\documentclass[12pt]{article}
\usepackage{amsmath,amssymb,amsthm,amsfonts,amscd,latexsym,verbatim}

\usepackage{indentfirst}
\usepackage{geometry}
\usepackage{microtype}

\usepackage{hyperref}

\usepackage[nottoc]{tocbibind}

\let\le\leqslant
\let\ge\geqslant

\let\mc\mathcal

\let\lf\lfloor
\let\rf\rfloor

\newcommand{\diff}{\mathrm{d}}
\newcommand{\NN}{{\mathbb Z_{>0}}}
\newcommand{\NNN}{{\mathbb Z_{\ge0}}}
\newcommand{\CC}{\mathbb C}

\newcommand{\ZZ}{\mathbb Z}
\newcommand{\QQ}{\mathbb Q}
\newcommand{\KK}{\mathbb K}
\newcommand{\ZK}{{\ZZ_\KK}}
\newcommand{\eps}{\varepsilon}
\renewcommand{\phi}{\varphi}
\newcommand{\ka}{\varkappa}
\renewcommand{\kappa}{\varkappa}

\DeclareMathOperator{\ord}{ord}

\newcommand{\house}[1]{\hbox{\vrule width.4pt%
 \vbox{{\hrule height.4pt}\vskip3pt\hbox{\,$\mathstrut #1$\,}}\vrule width.4pt}}

\newcommand{\MK}{{\mc M_\KK}}

\theoremstyle{plain}
\newtheorem{theorem}{Theorem}
\newtheorem{lemma}{Lemma}

\newtheorem{proposition}{Proposition}
\theoremstyle{remark}
\newtheorem{remark}{Remark}

\theoremstyle{definition}

\numberwithin{equation}{section}

\begin{document}

\hypersetup{pdfauthor={igor rochev},%
pdftitle={Linear independence measures for values of certain q-series}}

\title{Linear independence measures for values of certain $q$-series}

\author{I.~Rochev\thanks{Research is supported by RFBR (grant~No.~09-01-00371a).}}

\date{}

\maketitle

\begin{abstract}
We prove, in a quantitative form,  linear independence results for values of a certain class of $q$-series, which generalize classical $q$-hypergeometric series. These results refine our recent estimates.
\end{abstract}

\section{Main result}

Let $\KK$ be an algebraic number field of degree~$\ka=[\KK:\QQ]$, $\MK$ the set of all places of~$\KK$. For $v\in\MK$ we normalize the absolute value~${|\cdot|_v}$ as follows:
\begin{enumerate}
\item $|p|_v=p^{-1}$ for finite~$v|p$;
\item $|x|_v=|x|$ for~$x\in\QQ$ if $v|\infty$.
\end{enumerate}
Then for~$\alpha\in\KK^*$ we have the so-called product formula
\begin{equation*}
\prod_{v\in\MK}|\alpha|_v^{\ka_v}=1,
\end{equation*}
where $\ka_v=[\KK_v:\QQ_v]$ are the local degrees.

For a vector~$\vec\alpha=(\alpha_0,\alpha_1,\ldots,\alpha_n)\in\KK^{n+1}$ its (projective) absolute (multiplicative) height~$H(\vec\alpha)$ is given by
\begin{equation*}
H(\vec\alpha)=\prod_{v\in\MK}|\vec\alpha|_v^{\ka_v/\ka},\qquad|\vec\alpha|_v=\max\{|\alpha_0|_v,|\alpha_1|_v,\ldots,|\alpha_n|_v\}.
\end{equation*}
(In fact, by the product formula, $H(\cdot)$ is well defined on the projective space~$\mathbb{KP}^n$.) In particular, for~$\alpha\in\KK$ its absolute height is given by
\begin{equation*}
H(\alpha)=H\bigl((1,\alpha)\bigr)=\prod_{v\in\MK}\max\left\{|\alpha|_v^{\ka_v/\ka},1\right\}.
\end{equation*}
In view of the product formula, for~$\alpha\in\KK^*$ and any~$v\in\MK$ we have the so-called fundamental inequality
\begin{equation*}
\bigl|\log|\alpha|_v\bigr|\le\frac\ka{\ka_v}\log H(\alpha).
\end{equation*}

Suppose $q\in\KK$ and $w\in\MK$ satisfy $|q|_w>1$ and $|q|_v\le1$ for all~$v\in\MK\setminus\{w\}$. Further assume that polynomials~$P(x,y)\in\KK[x,y]$, $Q(x)\in\KK[x]$ satisfy $d:=\deg_yP(x,y)\ge1$ and $P(n,q^n)Q(n)\ne0$ for all~$n\in\NN$. Put
\begin{equation}\label{eq-Pinz}
\Pi_n(z)=\prod_{k=1}^nP(k,z^k)/Q(k)\qquad(n\in\NNN)
\end{equation}
and consider the function
\begin{equation*}
f(z)=\sum_{n=0}^\infty\frac{z^n}{\Pi_n(q)},\qquad z\in\CC_w,
\end{equation*}
where $\CC_w$ is the completion of the algebraic closure of~$\KK_w$.

The function~$f(z)$ is entire. Indeed, let
\begin{equation*}
P(x,y)=\sum_{\nu=0}^dp_\nu(x)y^\nu;
\end{equation*}
then for all sufficiently large~$n\in\NN$ we have
\begin{equation*}
|p_d(n)|_w\ge H\bigl(p_d(n)\bigr)^{-\ka/\ka_w}\ge n^{-c},\qquad c=\mathrm{const}.
\end{equation*}
Hence for large~$n$ we have
\begin{equation*}
|P(n,q^n)|_w\ge\frac12n^{-c}|q|_w^{dn}
\end{equation*}
and the assertion follows.

In this note we prove the following theorem.

\begin{theorem}\label{theorem-1}
Assume that the polynomials~$P(x,y),Q(x)$ satisfy (at least) one of the following two conditions:
\begin{description}
\item[(a)] $p_d(x)$ does not depend on~$x$, or
\item[(b)] $p_0(x)$ and $Q(x)$ do not depend on~$x$.
\end{description}
Let $m\in\NN$, $d_0\in\ZZ_{\ge d}$. Suppose numbers~$\alpha_j\in\KK^*$ and~$s_{j,k}\in\NN$ ($1\le j\le m$, $0\le k<d_0$) satisfy the following three conditions:
\begin{description}
\item[(i)] $\alpha_i\alpha_j^{-1}\notin q^\ZZ$ for $i\ne j$,
\item[(ii)] $s_{j,k}\le\deg p_d(x)$ for $1\le j\le m$ and $d\le k<d_0$, and
\item[(iii)] if $\deg p_0(x)=\deg Q(x)$, then $\alpha_j\notin(a/b)q^\NN$ for all~$j$, where $a$ and $b$ are the leading coefficients of~$p_0(x)$ and~$Q(x)$, respectively.
\end{description}
Then the numbers
\begin{equation*}
1,f^{(\sigma)}(\alpha_jq^k)\qquad(1\le j\le m,0\le k<d_0,0\le\sigma<s_{j,k})
\end{equation*}
are linearly independent over~$\KK$. Moreover, there exist (effective) positive constants~$C_0=C_0(P,Q,q)$ and~$H_0=H_0(P,Q,q,m,d_0,\alpha_j,s_{j,k})$ such that for any vector~$\vec A=(A_0,A_{j,k,\sigma})\in\KK^{1+\sum_{j,k}s_{j,k}}\setminus\{\vec0\}$ we have
\begin{equation*}
\left|A_0+\sum_{j=1}^m\sum_{k=0}^{d_0-1}\sum_{\sigma=0}^{s_{j,k}-1}A_{j,k,\sigma}f^{(\sigma)}(\alpha_jq^k)\right|_w\ge|\vec A\,|_w\exp\left(-C_0m^{2/3}(\log H)^{4/3}\right),
\end{equation*}
where $H=\max\{H(\vec A\,),H_0\}$.
\end{theorem}

In the case $\KK=\QQ$, $|\cdot|_w=|\cdot|$, $Q(x)=1$ the qualitative part of Theorem~\ref{theorem-1} was essentially proved by B\'ezivin~\cite{Bezivin-1990}. Moreover, B\'ezivin's result implies that in this case the corresponding assertion is valid (with $d_0=1$ and slightly more restrictive conditions posed on~$\alpha_j$) even if the polynomial~$P(x,y)$ does not satisfy conditions~(a)--(b) of Theorem~\ref{theorem-1}.

Recently the author~\cite{Rochev-2011} proposed a quantitative variant of B\'ezivin's method; in particular, a weak version of Theorem~\ref{theorem-1} was proved, with the estimate of the form~$\exp\bigl(-C_0m(\log H)^2\bigr)$. A modification of this method was proposed in~\cite{Rochev-preprint} for the case when the polynomials~$P(x,y),Q(x)$ do not depend on~$x$. In this case a much stronger result than Theorem~\ref{theorem-1} is valid: the estimate for the linear form is polynomial in~$H$ and the conditions posed on~$q$ can be weakened. In~\cite{Rochev-preprint} for simplicity only the case~$\KK=\QQ$, ${|\cdot|_w=|\cdot|}$ was considered but extension to the general case is straightforward (cf., e.\,g., \cite{Sankilampi-Vaananen-2007}).

Note that Theorem~\ref{theorem-1} allows one to describe all linear dependences (over~$\KK$) among values of the function~$f(z)$ and its derivatives at points of the field~$\KK$ (if the number~$q$ and the polynomials~$P(x,y),Q(x)$ satisfy the aforementioned conditions). Indeed, the function~$f(z)$ satisfies the functional equation
\begin{equation}\label{eq-functional-equation}
P\left(z\frac{\diff}{\diff z},D_q\right)\bigl(f(z)\bigr)=P(0,1)+Q\left(z\frac{\diff}{\diff z}\right)\bigl(zf(z)\bigr),\qquad D_qf(z):=f(qz),
\end{equation}
therefore, for any~$\alpha\in\KK^*$ and $s\ge\deg p_d(x)$ the number~$f^{(s)}(\alpha)$ can be expressed as a linear combination of the numbers $1$, $f^{(\sigma)}(\alpha)$ with $0\le\sigma<\deg p_d(x)$, and $f^{(\sigma)}(\alpha q^{-\nu})$ with $1\le\nu\le d$ and $\sigma\ge0$. It follows that, given numbers~$\beta_1,\ldots,\beta_l\in\KK^*$ and $t\in\NN$, there exist~$\alpha_j$ and $s_{j,k}$ satisfying the conditions of Theorem~\ref{theorem-1} such that the numbers~$f^{(\tau)}(\beta_j)$ ($1\le j\le l$, $0\le\tau<t$) can be expressed as linear combinations of~$1,f^{(\sigma)}(\alpha_jq^k)$ ($1\le j\le m$, $0\le k<d_0$, $0\le\sigma<s_{j,k}$). Hence any relation
\begin{equation*}
B_0+\sum_{j=1}^l\sum_{\tau=0}^{t-1}B_{j,\tau}f^{(\tau)}(\beta_j)=0
\end{equation*}
can be rewritten in the form
\begin{equation*}
A_0+\sum_{j=1}^m\sum_{k=0}^{d_0-1}\sum_{\sigma=0}^{s_{j,k}-1}A_{j,k,\sigma}f^{(\sigma)}(\alpha_jq^k)=0,
\end{equation*}
where $A_0,A_{j,k,\sigma}$ are certain linear combinations of~$B_0,B_{j,\tau}$. It follows from Theorem~\ref{theorem-1} that the coefficients~$B_0,B_{j,\tau}$ must satisfy the system of linear equations~$A_0=A_{j,k,\sigma}=0$. In other words, all linear dependences among values of the function~$f(z)$ and its derivatives at points of the field~$\KK$ follow from the functional equation~\eqref{eq-functional-equation}.

Theorem~\ref{theorem-1} is a simple consequence of results of~\cite{Rochev-2011}. Since the case~$\deg_xP(x,y)=\deg Q(x)=0$ was considered in~\cite{Rochev-preprint}, in what follows we assume that $\deg_xP(x,y)+\deg Q(x)>0$. In Section~\ref{section-summary} we summarize the required results from~\cite{Rochev-2011}. In Section~\ref{sec-main-proposition} we use them to construct auxiliary linear forms. In the final section Theorem~\ref{theorem-1} is proved.

\section{Summary}\label{section-summary}

Let $m\in\NN$, $d_0\in\ZZ_{\ge d}$, $\alpha_j\in\KK^*$ ($1\le j\le m$), $s_{j,k}\in\NN$ ($1\le j\le m$, $0\le k<d_0$). By~$\vec x$ denote the vector of variables~$\vec x=(x_0,x_{j,k,\sigma})$, where $1\le j\le m$, $0\le k<d_0$, $0\le\sigma<s_{j,k}$. Furthermore, put
\begin{gather*}
s_j=\max_{0\le k<d_0}s_{j,k}\qquad(1\le j\le m),\\
S=s_1+\ldots+s_m.
\end{gather*}

Consider the polynomials
\begin{gather}
\label{eq-un}u_n=u_n(z,\vec x)=\sum_{j=1}^m\sum_{k=0}^{d_0-1}\sum_{\sigma=0}^{s_{j,k}-1}\sigma!\binom n\sigma(\alpha_jz^k)^{n-\sigma}x_{j,k,\sigma}\in\KK[z,\vec x],\\
\label{eq-vn}v_n=v_n(z,\vec x)=\Pi_n(z)\cdot\left(x_0+\sum_{k=0}^n\frac{u_k(z,\vec x)}{\Pi_k(z)}\right)\in\KK[z,\vec x],
\end{gather}
where $\Pi_n(z)$ is given by~\eqref{eq-Pinz}.

\subsection{First case}\label{subsection-first-case}

Suppose the polynomial~$P(x,y)$ satisfies condition~(a) of Theorem~\ref{theorem-1}, i.\,e.,
\begin{equation*}
P(x,y)=p_dy^d+\sum_{\nu=0}^{d-1}p_\nu(x)y^\nu,\qquad p_d\in\KK^*.
\end{equation*}
Put $h=\deg Q(x)$,
\begin{equation*}
g_1=\max\left\{\max_{1\le\nu\le d}\frac{\deg p_{d-\nu}(x)}\nu,\frac hd\right\}>0.
\end{equation*}

Further, let~$\mc B$ be the backward shift operator given by
\begin{equation*}
\mc B\bigl(\xi(n)\bigr)=\xi(n-1).
\end{equation*}
For $1\le j\le m$ and $k\in\ZZ$ introduce the difference operator
\begin{equation*}
\mc A_{k,j}=\mc I-\alpha_jq^k\mc B,
\end{equation*}
where $\mc I$ is the identity operator, $\mc I\bigl(\xi(n)\bigr)=\xi(n)$.

Finally, for $l\ge0$, $n\ge(S+mh)l+m\sum_{k=0}^{l-1}\lf g_1k\rf$ put
\begin{equation*}
v_{l,n}(\vec x)=\prod_{k=0}^{l-1}\prod_{j=1}^m\mc A_{d_0-1-d-k,j}^{s_j+h+\lf g_1k\rf}\bigl(v_n(q,\vec x)\bigr):=\left(\prod_{k=0}^{l-1}\prod_{j=1}^m\mc A_{d_0-1-d-k,j}^{s_j+h+\lf g_1k\rf}\right)\bigl(v_n(q,\vec x)\bigr)\in\KK[\vec x].
\end{equation*}

Then we have the following lemma.

\begin{lemma}\label{lemma-1}
Assume that $\vec\omega=(\omega_0,\vec\omega_1)=(\omega_0,\omega_{j,k,\sigma})\in\CC_w^{1+\sum_{j,k}s_{j,k}}$ satisfies
\begin{equation*}
\omega_0+\sum_{j=1}^m\sum_{k=0}^{d_0-1}\sum_{\sigma=0}^{s_{j,k}-1}\omega_{j,k,\sigma}f^{(\sigma)}(\alpha_jq^k)=0.
\end{equation*}
Then for all $l\ge0$, $n\ge(S+mh)l+m\sum_{k=0}^{l-1}\lf g_1k\rf$ we have the estimate
\begin{equation*}
|v_{l,n}(\vec\omega)|_w\le|\vec\omega_1|_w|q|_w^{-ln+mg_1l^3/6+c(n+1)},
\end{equation*}
where the constant~$c>0$ depends only on the polynomials~$P,Q$ and the numbers~$q,m,d_0,\alpha_j,s_{j,k}$.
\end{lemma}

\begin{proof}
See~\cite[Lemma~3.2]{Rochev-2011}.
\end{proof}

\begin{remark}
Note that an estimate of the form
\begin{equation*}
|v_{l,n}(\vec\omega)|_w\ll_l|\vec\omega_1|_w|q|_w^{-ln+cn}
\end{equation*}
is trivial. Indeed, we have
\begin{equation*}
v_n(q,\vec\omega)=-\sum_{j=n+1}^\infty\frac{u_j(q,\vec\omega)}{\prod_{k=n+1}^jP(k,q^k)/Q(k)}=\sum_{j=1}^m\sum_{k=0}^{d_0-1}\sum_{\sigma=0}^{s_{j,k}-1}\omega_{j,k,\sigma}A_{j,k,\sigma},
\end{equation*}
where every~$A_{j,k,\sigma}$ has an asymptotic expansion of the form
\begin{equation*}
A_{j,k,\sigma}\sim\sum_{i=0}^\infty P_i(n)(\alpha_jq^{k-d-i})^n\qquad\text{as $n\to\infty$},
\end{equation*}
with $P_i(z)\in\KK[z]$, $\deg P_i\le \sigma+h+g_1i$. The required estimate follows immediately.
\end{remark}

\subsection{Second case}\label{subsection-second-case}

In this subsection we assume that the polynomials~$p_0(x)=p_0$ and $Q(x)=1$ do not depend on~$x$. Put
\begin{gather*}
\eps_0=\begin{cases}
0&\text{if $p_0=0$},\\
1&\text{if $p_0\ne0$},
\end{cases}\\
g_2=\max_{1\le\nu\le d}\frac{\deg p_\nu(x)}\nu>0.
\end{gather*}

For $0\le j\le m$ and $k\in\ZZ$ introduce the difference operator
\begin{equation*}
\mc A_{k,j}=\begin{cases}
\mc I-p_0z^k\mc B&\text{if $j=0$},\\
\mc I-\alpha_jz^k\mc B&\text{if $1\le j\le m$},
\end{cases}
\end{equation*}
where $\mc I,\mc B$ are same as above. Note that if~$p_0=0$, then $\mc A_{k,0}=\mc I$.

For $l\ge0$, $n\ge(S+\eps_0)l+(m+\eps_0)\sum_{k=0}^{l-1}\lf g_2k\rf$ put
\begin{equation*}
v_{l,n}=v_{l,n}(z,\vec x)=\prod_{k=0}^{l-1}\prod_{j=0}^m\mc A_{k,j}^{s_j+\lf g_2k\rf}(v_n)\in\KK[\vec x][z],\qquad s_0:=1.
\end{equation*}

Define the ($z$-)\emph{order} of a formal Laurent series~$\xi(z)=\sum_{n\in\ZZ}a_nz^n\ne0$ as
\begin{equation*}
\ord_z\xi(z)=\min\{n\mid a_n\ne0\};
\end{equation*}
furthermore, put $\ord_z0=+\infty$.

Then the following assertion holds.

\begin{lemma}\label{lemma-2}
For all $l\ge0$, $n\ge(S+\eps_0)l+(m+\eps_0)\sum_{k=0}^{l-1}\lf g_2k\rf$ we have
\begin{equation*}
\ord_zv_{l,n}\ge ln-(m+\eps_0)g_2l^3/6-c(n+1),
\end{equation*}
where the constant~$c>0$ depends only on the polynomial~$P$ and the numbers~$m,d_0,s_{j,k}$.
\end{lemma}

\begin{proof}
See~\cite[Lemma~3.3]{Rochev-2011}.
\end{proof}

\subsection{Non-vanishing lemma}

For $n\ge1$ put
\begin{equation*}
V_n=V_n(z,\vec x)=\det(v_{i+j})_{i,j=0}^{n-1}\in\KK[z,\vec x].
\end{equation*}

Then we have the following non-vanishing lemma.

\begin{lemma}\label{lemma-non-vanishing}
Assume that the polynomials~$P,Q$ and the numbers~$\alpha_j,s_{j,k}$ satisfy the conditions of Theorem~\ref{theorem-1}, $\vec\omega\in\CC_w^{1+\sum_{j,k}s_{j,k}}\setminus\{\vec0\}$. Then for any $n_0\in\NN$ there is an integer~$n$ within the range $n_0\le n\le c_1n_0+c_0$ such that $V_n(q,\vec\omega)\ne0$,  where $c_1=c_1(P,Q),c_0=c_0(Q,m,d_0,s_{j,k})$ are certain positive constants.
\end{lemma}

\begin{proof}
See~\cite[Lemma~4.3]{Rochev-2011}.
\end{proof}

\section{Main proposition}\label{sec-main-proposition}

We begin with some notation. Suppose
\begin{equation*}
A(\vec z\,)=A(z_1,\ldots,z_n)=\sum_{\vec\nu}A_{\vec\nu}z_1^{\nu_1}\ldots z_n^{\nu_n}\in\KK[\vec z\,]
\end{equation*}
is a polynomial; then for~$v\in\MK$ we put
\begin{equation*}
|A|_v=\begin{cases}
\sum_{\vec\nu}|A_{\vec\nu}|_v&\text{if $v|\infty$},\\
\max_{\vec\nu}|A_{\vec\nu}|_v&\text{if $v\nmid\infty$}.
\end{cases}
\end{equation*}
Furthermore, put
\begin{gather*}
H(A)=\prod_{v\in\MK}|A|_v^{\ka_v/\ka},\\
H_w(A)=\prod_{v\in\MK\setminus\{w\}}|A|_v^{\ka_v/\ka}.
\end{gather*}

In the following proposition we assume the hypotheses of Theorem~\ref{theorem-1}; we also keep the notation from the previous section.

\begin{proposition}\label{proposition-main}
There exists a constant~$g_0=g_0(P,Q)\in\NN$ such that for any positive integers~$l,n$ with~$n\ge mg_0l^2+Sl$ there is a linear form~$L_{l,n}(\vec x)\in\KK[\vec x]$ satisfying the following three conditions:
\begin{enumerate}
\item\label{condition-1-proposition} For any~$\vec\omega=(\omega_0,\vec\omega_1)=(\omega_0,\omega_{j,k,\sigma})\in\CC_w^{1+\sum_{j,k}s_{j,k}}$ such that
\begin{equation*}
\omega_0+\sum_{j=1}^m\sum_{k=0}^{d_0-1}\sum_{\sigma=0}^{s_{j,k}-1}\omega_{j,k,\sigma}f^{(\sigma)}(\alpha_jq^k)=0
\end{equation*}
we have
\begin{equation*}
|L_{l,n}(\vec\omega)|_w\le|\vec\omega_1|_w|q|_w^{-ln+mg_0l^3/2+O(n)}.
\end{equation*}
\item\label{condition-2-proposition} The following estimates are valid:
\begin{gather*}
H(L_{l,n})\le H(q)^{dn^2/2+O(n^{3/2})},\\
H_w(L_{l,n})\le H(q)^{O(n\log n)}.
\end{gather*}
\item\label{condition-3-proposition} For any~$\vec\omega\in\CC_w^{1+\sum_{j,k}s_{j,k}}\setminus\{\vec0\}$ and $l_0,n_0\in\NN$ with $n_0\ge mg_0l_0^2+Sl_0$ there exists an integer~$n$ with $n_0\le n\le g_0n_0+O(1)$ such that $L_{l_0,n}(\vec\omega)\ne0$.
\end{enumerate}
The constants in the Landau symbols~$O(\cdot)$ depend only on~$P,Q,q,m,d_0,\alpha_j,s_{j,k}$.
\end{proposition}

\begin{remark}
In fact, the inequality $n\le g_0n_0+O(1)$ in condition~\ref{condition-3-proposition} of Proposition~\ref{proposition-main} can be replaced by~$n\le n_0+O(l_0)$ (cf., e.\,g., \cite[Section~3]{Rochev-preprint}).
\end{remark}

The proof of the proposition is divided into two parts according to whether the polynomials~$P,Q$ satisfy condition~(a) or~(b) of Theorem~\ref{theorem-1}. Before we proceed let us make some preliminary remarks.

Without loss of generality we can assume that $Q(x)\in\ZK[x]$. Furthermore, let $D\in\NN$ be a general denominator of the numbers~$\alpha_j$ and coefficients of the polynomial~$P(x,y)$. Put
\begin{equation*}
I_n=D^n\prod_{k=1}^nQ(k)\in\ZK\setminus\{0\}.
\end{equation*}
It follows from~\eqref{eq-un}--\eqref{eq-vn} that $I_nv_n\in\ZK[z,\vec x]$. Moreover, $v_n$ is homogeneous in~$\vec x$ with~$\deg_{\vec x}v_n=1$, $\deg_zv_n=dn^2/2+O(n)$.

\subsection{Case~(a)}

We keep the notation from Subsection~\ref{subsection-first-case}. Let us show that we can take
\begin{equation*}
L_{l,n}(\vec x)=q^{\sum_{k=0}^{l-1}\sum_{j=1}^m(s_j+h+\lf g_1k\rf)(k+1)}v_{l,n}(\vec x).
\end{equation*}

It follows from Lemma~\ref{lemma-1} that condition~\ref{condition-1-proposition} of the proposition holds (provided that $g_0$ is large enough).

Further, since
\begin{equation*}
L_{l,n}(\vec x)=\prod_{k=0}^{l-1}\prod_{j=1}^m(q^{k+1}\mc I-\alpha_jq^{d_0-d}\mc B)^{s_j+h+\lf g_1k\rf}\bigl(v_n(q,\vec x)\bigr),
\end{equation*}
it is readily seen that $I_nL_{l,n}(\vec x)=A(q,\vec x)$ for some polynomial~$A(z,\vec x)\in\ZK[z,\vec x]$ with $\deg_zA\le dn^2/2+mg_1l^3/3+O(n)=dn^2/2+O(n^{3/2})$. Hence for finite~$v\in\MK$ we have
\begin{equation*}
|L_{l,n}|_v\le|I_n|_v^{-1}\max\{|q|_v,1\}^{dn^2/2+O(n^{3/2})}.
\end{equation*}

For a polynomial~$A(\vec z\,)=A(z_1,\ldots,z_n)=\sum_{\vec\nu}A_{\vec\nu}z_1^{\nu_1}\ldots z_n^{\nu_n}\in\KK[\vec z\,]$ put
\begin{equation*}
\mc L(A)=\sum_{\vec\nu}\house{A_{\vec\nu}},
\end{equation*}
where $\house{\alpha}=\max_{v|\infty}|\alpha|_v$ (in other words, $\house{\alpha}$ is the maximum of absolute values of $\alpha$'s conjugates). It follows from~\eqref{eq-un}--\eqref{eq-vn} that $\mc L(v_n)\le\exp\bigl(O(n\log n)\bigr)$. Hence for archimedean~$v\in\MK$ we have
\begin{equation*}
|v_n(q,\cdot)|_v\le\max\{|q|_v,1\}^{dn^2/2+O(n)}\exp\bigl(O(n\log n)\bigr).
\end{equation*}
This implies that
\begin{equation*}
|L_{l,n}|_v\le\max\{|q|_v,1\}^{dn^2/2+O(n^{3/2})}\exp\bigl(O(n\log n)\bigr).
\end{equation*}

Taking into account the estimate
\begin{equation*}
\prod_{v\nmid\infty}|I_n|_v^{-\ka_v/\ka}=\prod_{v|\infty}|I_n|_v^{\ka_v/\ka}\le\house{I_n}\le\exp\bigl(O(n\log n)\bigr)
\end{equation*}
and recalling that $|q|_v\le1$ for all $v\ne w$, we obtain condition~\ref{condition-2-proposition} of the proposition.

Condition~\ref{condition-3-proposition} follows from Lemma~\ref{lemma-non-vanishing}. Indeed, if $v_{l_0,n}(\vec\omega)=0$ for all~$n$ with $n_0\le n\le N_0$, then we have $V_n(q,\vec\omega)=0$ for $n_0+1\le n\le N_0/2+1$.

This concludes the proof of the proposition in the first case.

\subsection{Case~(b)}

Without loss of generality, we can assume that $Q(x)=1$.

We keep the notation from Subsection~\ref{subsection-second-case}. Let us show that we can take
\begin{equation*}
L_{l,n}(\vec x)=q^{-\ord_zv_{l,n}}v_{l,n}(q,\vec x).
\end{equation*}

Suppose $\vec\omega=(\omega_0,\vec\omega_1)=(\omega_0,\omega_{j,k,\sigma})\in\CC_w^{1+\sum_{j,k}s_{j,k}}$ satisfies
\begin{equation*}
\omega_0+\sum_{j=1}^m\sum_{k=0}^{d_0-1}\sum_{\sigma=0}^{s_{j,k}-1}\omega_{j,k,\sigma}f^{(\sigma)}(\alpha_jq^k)=0.
\end{equation*}
This implies that
\begin{equation*}
|v_n(q,\vec\omega)|_w=\left|-\sum_{j=n+1}^\infty\frac{u_j(q,\vec\omega)}{\prod_{k=n+1}^jP(k,q^k)/Q(k)}\right|_w\le|\vec\omega_1|_w|q|_w^{O(n)}.
\end{equation*}
Therefore,
\begin{equation*}
|v_{l,n}(q,\vec\omega)|_w\le|\vec\omega_1|_w|q|_w^{(m+\eps_0)g_2l^3/3+O(n)},
\end{equation*}
and condition~\ref{condition-1-proposition} of the proposition follows from Lemma~\ref{lemma-2}.

Further, we have
\begin{gather*}
I_nv_{l,n}\in\ZK[z,\vec x],\\
\deg_zv_{l,n}\le dn^2/2+O(n^{3/2}).
\end{gather*}
Therefore, if $v\nmid\infty$, then
\begin{equation*}
|L_{l,n}|_v\le|I_n|_v^{-1}\max\{|q|_v,1\}^{dn^2/2+O(n^{3/2})}.
\end{equation*}
Moreover,
\begin{equation*}
\mc L(v_{l,n})\le\exp\bigl(O(n\log n)\bigr),
\end{equation*}
hence for~$v|\infty$ we have
\begin{equation*}
|L_{l,n}|_v\le\max\{|q|_v,1\}^{dn^2/2+O(n^{3/2})}\exp\bigl(O(n\log n)\bigr).
\end{equation*}
These estimates imply condition~\ref{condition-2-proposition} of the proposition.

Condition~\ref{condition-3-proposition} follows from Lemma~\ref{lemma-non-vanishing}.

Proposition~\ref{proposition-main} is proved.

\section{Proof of Theorem~\ref{theorem-1}}

Take $n_0=mg_0l^2+Sl$, where $l\in\NN$ will be chosen later. It follows from~Proposition~\ref{proposition-main} that there is an integer~$n$ with $mg_0l^2+O(l)\le n\le mg_0^2l^2+O(l)$ such that $L_{l,n}(\vec A\,)\ne0$. By the product formula, we have
\begin{equation*}
|L_{l,n}(\vec A\,)|_w^{\ka_w/\ka}=\prod_{v\in\MK\setminus\{w\}}|L_{l,n}(\vec A\,)|_v^{-\ka_v/\ka}\ge\prod_{v\in\MK\setminus\{w\}}\bigl(|L_{l,n}|_v|\vec A\,|_v\bigr)^{-\ka_v/\ka}=|\vec A\,|_w^{\ka_w/\ka}\bigl(H_w(L_{l,n})H(\vec A\,)\bigr)^{-1},
\end{equation*}
therefore,
\begin{equation*}
|L_{l,n}(\vec A\,)|_w\ge|\vec A\,|_w\bigl(H_w(L_{l,n})H\bigr)^{-\ka/\ka_w}\ge|\vec A\,|_wH^{-\ka/\ka_w}|q|_w^{O(l^2\log l)}.
\end{equation*}

Let $\vec\omega=(\omega_0,\omega_{j,k,\sigma})$ be given by
\begin{gather*}
\omega_{j,k,\sigma}=A_{j,k,\sigma},\\
\omega_0=-\sum_{j=1}^m\sum_{k=0}^{d_0-1}\sum_{\sigma=0}^{s_{j,k}-1}\omega_{j,k,\sigma}f^{(\sigma)}(\alpha_jq^k).
\end{gather*}
It follows from Proposition~\ref{proposition-main} that
\begin{equation*}
|L_{l,n}(\vec\omega)|_w\le|\vec A\,|_w|q|_w^{-ln+mg_0l^3/2+O(n)}\le|\vec A\,|_w|q|_w^{-mg_0l^3/2+O(l^2)}.
\end{equation*}
Take $l=\Bigl\lf\left(\frac{3\log H}{mg_0\log H(q)}\right)^{1/3}\Bigr\rf$. Then we have
\begin{equation*}
|L_{l,n}(\vec\omega)|_w\le\frac12|L_{l,n}(\vec A\,)|_w,
\end{equation*}
provided $H_0$ (and hence $H$) is sufficiently large. Therefore,
\begin{equation*}
|L_{l,n}(\vec A\,)-L_{l,n}(\vec\omega)|_w\ge\frac12|L_{l,n}(\vec A\,)|_w\ge\frac12|\vec A\,|_w\bigl(H_w(L_{l,n})H\bigr)^{-\ka/\ka_w}.
\end{equation*}

On the other hand, we have
\begin{equation*}
|L_{l,n}(\vec A\,)-L_{l,n}(\vec\omega)|_w\le|L_{l,n}|_w|A_0-\omega_0|_w.
\end{equation*}
Thus,
\begin{equation*}
\left|A_0+\sum_{j=1}^m\sum_{k=0}^{d_0-1}\sum_{\sigma=0}^{s_{j,k}-1}A_{j,k,\sigma}f^{(\sigma)}(\alpha_jq^k)\right|_w=|A_0-\omega_0|_w\ge\frac12|\vec A\,|_w\bigl(H(L_{l,n})H\bigr)^{-\ka/\ka_w},
\end{equation*}
and Theorem~\ref{theorem-1} follows from Proposition~\ref{proposition-main}.

\end{document}